\documentclass[11pt,titlepage,a4paper]{article}
\usepackage{amsfonts}
\usepackage{latexsym,emlines}

\newcommand{\QED}[1][6]{\hbox{\vrule width #1pt height #1pt depth 0pt}}

\newtheorem{Th}{Theorem}[section]		
\newtheorem{Prop}{Proposition}[section]		
\newtheorem{Lem}{Lemma}[section]		
\newtheorem{Def}{Definition}[section]		

\newtheorem{exam}{Example}[section]		
\newenvironment{Exam}[1]{\begin{exam}\label{#1}\slshape}{\end{exam}}

\input amssym.def \newsymbol\varkappa 207B
\input{amssym}

\begin{document}

\title{\bfseries\vspace*{-1.8in}
\vspace*{-7ex}
{
\begin{flushright}
	  \textbf{\large LANL xxx E-print archive No. dg-ga/9709016}\\[5ex]
\end{flushright}
}
\bfseries\huge Transports along maps \\[6pt] in fibre bundles
}
\vspace{10pt}
\author{Bozhidar Z. Iliev
\thanks{Permanent address:
Department of Mathematical Modeling,
Institute for Nuclear Research and \mbox{Nuclear} Energy,
Bulgarian Academy of Sciences,
Boul. Tzarigradsko chauss\'ee~72, 1784 Sofia, Bulgaria}
\thanks{E-mail address: bozho@inrne.bas.bg}
\thanks{URL: http://www.inrne.bas.bg/mathmod/bozhome/}
}

\newlength{\bo}\newlength{\ho}\newlength{\up}
\newlength{\down}\newlength{\middle}
\newcommand{\BOZHO}{\leavevmode\hbox{\slshape\bfseries
\settowidth{\bo}{BO}\settowidth{\ho}{HO}%
\settowidth{\middle}{/}\settoheight{\up}{BOZHO}%
\settodepth{\down}{/}%
\addtolength{\up}{+0.15\up}\addtolength{\bo}{+\middle}%
\rule[\up]{\bo}{0.15ex}\hspace{-\bo}BO%
\hspace{+0.09em}\raisebox{+0.17\up}{/}%
\hspace{-0.20em}\raisebox{+0.71\up}{$\bullet$}%
\hspace{-0.33em}\hspace{-1.14\middle}\raisebox{-0.4\up}{$\bullet$}%
\hspace{-0.30em}%
\addtolength{\down}{-0.41\down}\addtolength{\ho}{+1.5\middle}%
\rule[-\down]{\ho}{0.15ex}\addtolength{\ho}{-\middle}%
\hspace{-\ho}\hspace{+0.18em}\raisebox{+0.17\up}{HO}%
}}


\date{
 Ended: December 16, 1995 \\
 Revised: November 11, December 23, 1996\\[4pt]
 Updated: July 4, 1998\\
 Produced: \fbox{\today}\\[1ex]
 Submitted to JINR communication: January 13, 1997\\
Published: JINR Communication E5-97-2, 1997\\[1.4 ex]
LANL xxx archive server E-print No. dg-ga/9709016\\[7ex]
{\Huge\BOZHO}
}

\renewcommand{\thefootnote}{\fnsymbol{footnote}}
\maketitle			
\renewcommand{\thefootnote}{\arabic{footnote}}

\tableofcontents

   \begin{abstract}

   	Generalized are the investigated in other works of the author
transports along paths in fibre bundles to transports along arbitrary maps in
them. Their structure and some properties are studied. Special attention
is paid on the linear case and on the one when map's domain is a Cartesian
product of two sets. There are considered the consistency with bundle
morphisms and a number of special cases.

   \end{abstract}

\section {\bfseries Introduction}
\label{I}
\setcounter {equation} {0}

	In previous papers
(see, e.g.,~\cite{f-TP-general,f-LTP-general}) we
have studied the transports along paths in fibre bundles. In them is
not always essential the fact that the transports are along paths.
This suggests a way of generalizing these investigations which
is the subject of the present work.

	Sect.~\ref{II} gives and discusses the basic definition of
transports along maps in fibre bundles. Sect.~\ref{III} studies in
details the case when the map's domain is a Cartesian product of two
sets. Here presented are certain examples too. Sect.~\ref{IV} is
devoted to linear transports along maps in vector bundles.
Partial derivations along maps are introduced as well as
the general concepts of curvature and
torsion. It is pointed out how
a number of the already obtained results
concerning linear transports along paths.
can \emph{mutatis mutandis} be transferred in the
investigated here general case
Sect.~\ref{V} investigated, in analogy
with~\cite{f-TP-morphisms}, the consistency
(compatibility) of transports along maps in fibre bundles with bundle
morphisms between them.  Sect.~\ref{VI} closes the paper with a
discussion of different problems: An interpretation is given of the
obtained in Sect.~\ref{III} decomposition of transports along maps
whose domain is a Cartesian product of two sets. A scheme is proposed
for performing operations with elements of different fibres of a bundle
as well as with its
sections. It is proved that the Hermitian metrics on a differentiable
manifold are in one to one correspondence with the transports along the
identity map in an appropriate tensor bundle over it. At the end, some
remarks concerning tensor densities are discussed.

\section {\bfseries The basic definition.\\
				Special cases and discussion}
\label{II}
\setcounter {equation} {0}

	The fact that $\gamma$ is a path in definition~2.1
of~\cite{f-TP-general} for a transport along paths in fibre bundles is
insignificant from a logical view-point. This observation, as well as
other reasons, leads to the following generalization.

	Let $(E,\pi,B)$ be a topological fibre bundle with base $B$,
total space $E$, projection $\pi:E\to B$, and homeomorphic fibres
$\pi^{-1}(x),\ x\in B$.
Let the set $N$ be not empty ($N\neq\emptyset$) and there  be
given a map $\varkappa:N\to B$. By $id_M$ is denoted the identity map
of the set $M$.

	\begin{Def} \label{Def2.1}
	A transport along maps in the fibre bundle $(E,\pi,B)$ is a map
$K$ assigning to any map $\varkappa:N\to B$ a map  $K^\varkappa$,
transport along $\varkappa$, such that
$K^\varkappa:(l,m)\mapsto K_{l\to m}^{\varkappa}$, where for
every $l,m\in N$ the map
\begin{equation}	\label{2.1}
K_{l\to m}^{\varkappa}:\pi^{-1}(\varkappa(l))\to \pi^{-1}(\varkappa(m)),
\end{equation}
called transport along  $\varkappa$ from  $l$ to $m$, satisfies the
equalities:
\begin{eqnarray}
K_{m\to n}^{\varkappa} \circ  K_{l\to m}^{\varkappa} &=&
K_{l\to n}^{\varkappa} , \quad l,m,n\in N,  		\label{2.2}\\
K_{l\to l}^{\varkappa} &=& id_{\pi^{-1}(\varkappa(l))},
\quad l\in N.    					\label{2.3}
\end{eqnarray}
\end{Def}

	The formal analogy of this definition with definition~2.1
of~\cite{f-TP-general} is evident. In particular, if $\varkappa$ is a
path in $B$, i.e. if $N$ is an  $\mathbb{R}$-interval, the former
definition reduces to the latter.
The two definitions coincide also in the `flat' case when
$N=B$ and  $\varkappa=id_B$. In fact, in this case
 $I_{s\to t}^{\gamma}:=K_{\gamma(s)\to\gamma(t)}^{id_B}$
for a path $\gamma:J\to B$, $J$ being an
$\mathbb{R}$-interval, $s,t\in J$,
defines a transport along paths in $(E,\pi,B)$ which depends only on
the points $\gamma(s)$ and  $\gamma(t)$ but not on the path
$\gamma$ itself. On the opposite, if
$I$ is a transport along paths having the last property, then
$K_{\gamma(s)\to\gamma(t)}^{id_B}:=I_{s\to t}^{\gamma}$
is a transport along the identity map of $B$ in $(E,\pi,B)$.
By~\cite[theorem~6.1]{f-LTP-Cur+Tor-prop} the so defined transports
along paths are flat, i.e. their curvature vanishes in the case when
they are linear and $(E,\pi,B)$ is a vector bundle. Due to these facts,
we call the transports along the identity map \emph{flat transports}.

	The general form of a transport along maps is given by
	\begin{Th}	\label{Th2.1}
Let $\varkappa:N\to B$. The map
$K:\varkappa\mapsto K^\varkappa:(l,m)\mapsto K_{l\to m}^{\varkappa}$,
$l,m\in N$
is a transport along $\varkappa$ if and only if there exist a set
$Q$ and a family of one-to-one maps
 $\{ F_{n}^{\varkappa}:\pi^{-1}(\varkappa(n))\to Q,\ n\in N \}$
such that
	\begin{equation}	\label{2.4}
K_{l\to m}^{\varkappa} = \left(F_{m}^{\varkappa}\right)^{-1} \circ
\left(F_{l}^{\varkappa}\right),	\quad l,m\in N.
	\end{equation}
The maps $F_{n}^{\varkappa}$ are defined up to a left composition with
1:1 map depending only on $\varkappa$, i.e.~(\ref{2.4}) holds for given
families of maps
$\{ F_{n}^{\varkappa}:\pi^{-1}(\varkappa(n))\to Q,\ n\in N \}$ and
\(
\{ ^\prime\! F_{n}^{\varkappa}:\pi^{-1}(\varkappa(n))\to
^\prime\! Q,\ n\in N \}
\)

for some sets $Q$ and  $^\prime\! Q$ iff there is 1:1 map
$D^\varkappa:Q\to^\prime\! Q$ such that
	\begin{equation}	\label{2.5}
^\prime\! F_{n}^{\varkappa} = D^\varkappa \circ F_n^\varkappa,\quad n\in N.
	\end{equation}
	\end{Th}

	\emph{Proof.} This theorem is a trivial corollary of lemmas~3.1
and~3.2 of~\cite{f-TP-general} for $Q_n=\pi^{-1}(\varkappa(n)),\ n\in N$
 and  $R_{l\to m}=K_{l\to m}^{\varkappa},\ l,m\in N.\>$\QED

	The formal analogy is evident between transports along maps and
the ones along paths. The causes for this are definition~\ref{Def2.1}
of  this work and~\cite[definition 2.1]{f-TP-general}, as
 well as theorem~\ref{Th2.1} of the present work
and~\cite[theorem~3.1]{f-TP-general}.
Due to this almost all results concerning
transports along paths are valid \emph{mutatis mutandis} for transports
along maps.  Exceptions are the results which use explicitly the fact
that a path is a map from a real interval to a certain set, viz.  in
which special properties of the
$\mathbb{R}$-intervals, such as ordering, the Abelian structure (the
operation addition) an so on, are used.  This transferring of results
can formally be done by substituting the symbols $\varkappa$ for
$\gamma$, $N$ for $J$, $K$ for $I$, $l,m,n\in N$ for $r,s,t\in J$,
and the word map(s) for the word path(s).

	For example, definition~2.2, proposition~2.1 and example~2.1
of~\cite{f-TP-general} now read:

\begin{Def} 	\label{Def2.2}
	The section  $\sigma\in\mathrm{Sec}(E,\pi,B)$ undergoes a
(K-)transport or is (K-)transported (resp. along  $\varkappa:N\to B$)
if the equality
\begin{equation}					\label{2.6}
\sigma(\varkappa(m)) =
K_{l\to m}^{\varkappa}\sigma(\varkappa(l)),\quad l,m\in N
\end{equation}
holds for every (resp. the given) map $\varkappa:N\to B$.
\end{Def}

\begin{Prop}	\label{Prop2.1}
	If~(\ref{2.6}) holds for a fixed $l\in N$, then
it is valid for every $l\in N$.
\end{Prop}

		\begin{Exam}{Exam2.1}
		If $(E,\pi,B)$ has a foliation structure $\{K_\alpha;\
\alpha\in A\}$, then the lifting $\overline{\varkappa}_u:N\to E$ of
$\varkappa:N\to B$  through $u\in E$ given by
\[
\overline{\varkappa}_u(l) :=
		\pi^{-1}(\varkappa(l)) \bigcap K_{\alpha(u)},
\]
where $K_{\alpha(u)}\ni u$, defines a transport  $K$ along maps through
\begin{displaymath}
K_{l\to m}^{\varkappa}(u) := \overline{\varkappa}_u(m),\quad
u\in\pi^{-1}(\varkappa(l)),\quad l,m\in N.
\end{displaymath}
		\end{Exam}

	On the transports along maps additional restrictions can be
imposed, such as (cf.~\cite[Sect.~2.2 and~2.3]{f-TP-general}):

	\begin{itemize}

\item
	the locality condition:
	\begin{equation}	\label{2.7}
K_{l\to m}^{\varkappa|N^\prime} = K_{l\to m}^{\varkappa},\quad
l,m\in N^\prime\subset N;
	\end{equation}

\item
	the `reparametrization invariance' condition:
	\begin{equation}	\label{2.8}
K_{l\to m}^{\varkappa \circ \tau} =
			K_{\tau(l)\to \tau(m)}^{\varkappa},\quad
\mathrm{\ for\  1:1\  map}\ \tau:N''\to N,\quad l,m\in N'';
	\end{equation}

\item
	the consistency with a bundle binary operation
 $\beta: x\to\beta_x,\ x\in B$ (e.g. a metric, i.e. a scalar product):
	\begin{equation}	\label{2.9}
\beta_{\varkappa(l)} =
		\beta_{\varkappa(m)}\circ\left(K_{l\to m}^{\varkappa}
\times K_{l\to m}^{\varkappa}\right)
	\end{equation}
for \( \beta_{x}:\pi^{-1}(x)\times\pi^{-1}(x)\to M,\ x\in B,\)
$M$ being a set, e.g. $M=\mathbb{R,\ C}$;

\item
	the consistency with the vector structure of a complex (or real)
vector bundle:
	\begin{equation}      \label{2.10} \!\! \!\! \!\! \!\! \!\! \!\!
K_{l\to m}^{\varkappa}(\lambda u + \mu v) =
\lambda K_{l\to m}^{\varkappa}u + \mu K_{l\to m}^{\varkappa}v,\ \>
\lambda,\mu\in\mathbb{C\mathrm{\ (or\ }R)},\ \>
u,v\in\pi^{-1}(\varkappa(l)).
	\end{equation}

\end{itemize}
	The last condition defines the set of
\emph{linear transports along maps}.

	Examples of results that do not have analogs in our general
case are propositions~2.3, 3.3, and~3.4 of~\cite{f-TP-general}.
But the all definitions and results of Sect.~4 and Sect.~5
of~\cite{f-TP-general} have analogs in this case.
They can be obtained by making the above-pointed substitutions.

\section {\bfseries The composite case}
\label{III}
\setcounter {equation} {0}

	Of special interest are transports along maps whose domain
has a structure of a Cartesian product, i.e. maps like
$\varkappa:N\to B$
with $N=A\times M$, $A$ and $M$ being not empty sets.
In this section, transports
$K_{(a,x)\to(b,y)}^{\eta}$ along $\eta:A\times M\to B$ are considered
and their general form is found. Here and below $a,b,c\in A$ and
$x,y,z\in M$.
By $\eta(\cdot,x):A\to B$ and $\eta(a,\cdot):M\to B$ are denoted,
respectively, the maps $\eta(\cdot,x):a\mapsto\eta(a,x)$ and
$\eta(a,\cdot):x\mapsto\eta(a,x)$.

	Applying~(\ref{2.2}), we get
	\begin{equation} 	\label{3.1}
K_{(a,x)\to(b,y)}^{\eta} =
K_{(a,y)\to(b,y)}^{\eta} \circ K_{(a,x)\to(a,y)}^{\eta} =
K_{(b,x)\to(b,y)}^{\eta} \circ K_{(a,x)\to(b,x)}^{\eta}.
	\end{equation}

	Using~(\ref{2.1})--(\ref{2.3}), we see that
 $^x\! K_{a\to b}^{\eta} := K_{(a,x)\to(b,x)}^{\eta} $ and
 $_a K_{x\to y}^{\eta} := K_{(a,x)\to(a,y)}^{\eta} $
satisfy~(\ref{2.1})--(\ref{2.3}) with, respectively, $\varkappa=\eta$,
$K={^x\!}K$, $l=a$, $m=b$, and $n=c$ and $\varkappa=\eta$, $K={_a}K$,
$l=x$, $m=y$, and $n=z$.  Consequently, a transport along $\eta$
decomposes to a composition of two (commuting) maps
satisfying~(\ref{2.1})--(\ref{2.3}).
Note that if the locality condition~(\ref{2.7}) holds, then these maps
are simply the transports along $\eta(\cdot,x)$ and $\eta(a,\cdot)$.

	So, applying lemma 3.1 of~\cite{f-TP-general}, we find
	\begin{equation}	\label{3.2}
K_{(a,x)\to(b,x)}^{\eta} =
\left(
^x\! H_{b}^{\eta}\right)^{-1} \circ \left(^x\! H_{a}^{\eta}
\right),
\qquad K_{(a,x)\to(a,y)}^{\eta} =
\left(
^a\! G_{y}^{\eta}\right)^{-1} \circ \left(^a\! G_{x}^{\eta}
\right),
	\end{equation}
where $^x\! H_{a}^{\eta}:\pi^{-1}(\eta(a,x))\to Q_H$ and
$^a\! G_{x}^{\eta}:\pi^{-1}(\eta(a,x))\to Q_G$
are 1:1 maps on some sets $Q_H$ and $Q_G$ respectively.
(The maps $^x\! H_{a}^{\eta}$ and
$^a\! G_{x}^{\eta}$ are defined up to a left composition
with 1:1 maps depending on the pairs $x$ and $\eta$ and  $a$ and $\eta$
respectively - see~\cite[lemma 3.2]{f-TP-general}.)

	The substitution of~(\ref{3.2}) into~(\ref{3.1}) yields
	\begin{equation}   \label{3.3}
	\begin{array}{rcl}
K_{(a,x)\to(b,y)}^{\eta} &=&
\left(^y\!H_{b}^{\eta}\right)^{-1}\circ \left(^y\!H_{a}^{\eta}\right)
\circ
\left(^a\!G_{y}^{\eta}\right)^{-1}\circ \left(^a\!G_{x}^{\eta}\right) =
\\  &=&
\left(^b\!G_{y}^{\eta}\right)^{-1}\circ \left(^b\!G_{x}^{\eta}\right)
\circ
\left(^x\!H_{b}^{\eta}\right)^{-1}\circ \left(^x\!H_{a}^{\eta}\right).
	\end{array}
	\end{equation}

Separating the terms depending on $x$ and  $y$ in the second equality,
we see that there exist one-to-one maps
$C_{a\to b}^{\eta}:Q_G\to Q_G$ which are independent of $x$
and such that
	\begin{equation}	\label{3.4}
\left(^b\!G_{x}^{\eta}\right) \circ \left(^x\!H_{b}^{\eta}\right)^{-1}
\circ
\left(^x\!H_{a}^{\eta}\right) \circ \left(^a\!G_{x}^{\eta}\right)^{-1} =
C_{a\to b}^{\eta} .
	\end{equation}

	It is trivial to check the equalities
 $C_{a\to b}^{\eta}=C_{c\to b}^{\eta} \circ C_{a\to c}^{\eta}$ and
 $C_{a\to a}^{\eta} = id_{Q_G}$.
Hence, by~\cite[lemma 3.1]{f-TP-general}, we have
 $C_{a\to b}^{\eta} = \left(C_{b}^{\eta}\right)^{-1} \circ C_{a}^{\eta}$
with certain 1:1 maps $C_{a}^{\eta}:Q_G\to Q_C$
(defined up to a left composition with a map depending only on
$\eta$~\cite[lemma 3.2]{f-TP-general}) on some set $Q_C$ .
The substitution of this result into~(\ref{3.4}) and the separation of
the terms depending on $a$ and $b$, shows the existence of 1:1 map
 $D_{x}^{\eta}:Q_H\to Q_C$ depending on  $\eta$ and $x$,
for which
	\begin{equation}	\label{3.5}
\left(^x\!H_{a}^{\eta}\right) \circ \left(^a\!G_{x}^{\eta}\right)^{-1}
\circ
\left(C_{a}^{\eta}\right)^{-1} = \left(D_{x}^{\eta}\right)^{-1} .
	\end{equation}
Hereout
	\begin{equation}	\label{3.6}
^x\!H_{a}^{\eta} = \left(D_{x}^{\eta}\right)^{-1} \circ
 C_{a}^{\eta} \circ \left(^a\!G_{x}^{\eta}\right)
\quad\mathrm{or}\quad
^a\!G_{x}^{\eta} = \left(C_{a}^{\eta}\right)^{-1}  \circ
D_{x}^{\eta} \circ \left(^x\!H_{a}^{\eta}\right) .
	\end{equation}
Substituting~(\ref{3.6}) into~(\ref{3.3}), we finally, in accordance
with~(\ref{2.4}), get
	\begin{equation}	\label{3.7}
K_{(a,x)\to(b,y)}^{\eta} =
\left(F_{(b,y)}^{\eta}\right)^{-1} \circ F_{(a,x)}^{\eta}
	\end{equation}
with
	\begin{equation}	\label{3.8}
F_{(a,x)}^{\eta} = C_{a}^{\eta} \circ \left(^a\!G_{x}^{\eta}\right) =
D_{x}^{\eta} \circ \left(^x\!H_{a}^{\eta}\right)
:\pi^{-1}(\eta(a,x))\to Q_C .
	\end{equation}

 	As we noted above, the maps
$^a\!G_{x}^{\eta},\ ^x\!H_{a}^{\eta}$, and
$C_{a}^{\eta}$ are defined up to the changes
	\begin{equation}	\label{3.9}
{^a\!}G_{x}^{\eta} \to
{^a\!}P_{G}^{\eta}\circ\left({^a\!}G_{x}^{\eta}\right),\qquad
^x\!H_{a}^{\eta} \to {^x\!}P_{H}^{\eta} \circ
\left({^x\!}H_{a}^{\eta}\right),
	\end{equation}
and
$ C_{a}^{\eta} \to P_{C}^{\eta} \circ C_{a}^{\eta} $, respectively,
where
\(
{^a\!}P_{G}^{\eta}
:Q_G \to Q_G,\ {^x\!}P_{H}^{\eta} :Q_H \to Q_H,\ P_{G}^{\eta}
						:Q_C \to Q_C
\)
are 1:1 mappings.
	The transformation concerning $C_{a}^{\eta}$ is valid if
$C_{a\to b}^{\eta}$ is defined independently. But this is not our case.
Due to~(\ref{3.4}) the changes~(\ref{3.9}) imply
$C_{a\to b}^{\eta} \to {^b\!}P_{G}^{\eta} \circ C_{a\to b}^{\eta}\circ
\left({^a\!}P_{G}^{\eta}\right)^{-1}$.
To describe this transformation through $C_{a}^{\eta}$ we must have
	\begin{equation}	\label{3.10}
C_{a}^{\eta} \to
P_{C}^{\eta} \circ C_{a}^{\eta} \circ
			\left({^a\!}P_{G}^{\eta}\right)^{-1} .
	\end{equation}
	From~(\ref{3.5}) it is easy to verify that the
transformations~(\ref{3.9}) and~(\ref{3.10}) imply
	\begin{equation}	\label{3.101}
D_{x}^{\eta} \to
P_{C}^{\eta} \circ D_{x}^{\eta} \circ
				\left(^x\!P_{H}^{\eta}\right)^{-1}.
	\end{equation}

     At the end, according to~(\ref{3.8}), all this leads to the change
	\begin{equation}	\label{3.11}
F_{(a,x)}^{\eta} \to P_{C}^{\eta} \circ F_{(a,x)}^{\eta} ,
	\end{equation}
as should be by~(\ref{2.5}).

	Together~(\ref{3.9})--(\ref{3.11}) form the set of
transformations under which our theory is invariant.

	Thus we have proved
	\begin{Prop}	\label{Prop3.1}
	The set of maps $\{K_{(a,x)\to(b,y)}^{\eta}\}$ forms a transport
along $\eta:A\times M\to B$ iff~(\ref{3.7}) and~(\ref{3.8})
are valid for some 1:1 maps shown on the commutative diagram
	\begin{equation}	\label{3.12}
\begin{array}{c}
\input{alongmap.pic}
\end{array}
	\end{equation}
that are defined up to the transformations given
by~(\ref{3.9})--(\ref{3.11}).
	\end{Prop}

	\textbf{Remark.}
In fact  $^x\!H_{a}^{\eta}$ and $^a\!G_{x}^{\eta}$
determine the `restricted' transports $K_{(a,x)\to(b,x)}^{\eta}$ and
 $K_{(a,x)\to(a,y)}^{\eta}$ through~(\ref{3.2}). In the case when the
locality condition~(\ref{2.7}) holds for $\varkappa=\eta$, they are
equal, respectively, to the transports $K_{a\to b}^{\eta(\cdot,x)}$
and $K_{x\to y}^{\eta(a,\cdot)}$
along the restricted maps
$\eta(\cdot,x)$ and $\eta(a,\cdot)$.
Note also that if $Q_G,\ Q_H,\ \mathrm{and}\ Q_C$ are
regarded as different  typical fibres of $(E,\pi,B)$, then the shown
maps represent different ways for mapping a concrete fibre on them.
This interpretation is more natural if, one puts
$Q_G= Q_H= Q_C=Q$, $Q$ being the typical fibre of $(E,\pi,B)$.
This is possible
due to the arbitrariness in
$^a\!G_{x}^{\eta},\ ^x\!H_{a}^{\eta}$ and $C_{a}^{\eta}$.

		\begin{Exam}{Exam3.1}
		Now we shall prove that the considered
in~\cite{f-Jacobi} transport
in a family of (vector) bundles
$\{\xi_a:\ \xi_a=(E_a,\pi_a,M),\ a\in A\}$
over one and the same base (manifold)  $M$ defined by the maps
 ${^{a,b}\!}I_{x\to y}:\pi_a^{-1}(x)\to \pi_b^{-1}(y)$, such that
 ${^{b,c}\!}I_{y\to z}\circ{^{a,b}\!}I_{x\to y} = {^{a,c}\!}I_{x\to z}$
and  ${^{a,a}\!}I_{x\to x} = id_{\pi_{a}^{-1}(x)}$,
is a (flat) transport along the identity map of the base of a suitably
chosen fibre bundle.

	A given family  $\{(E_a,\pi_a,M),\ a\in A\}$ of fibre bundles
over one and the same base is equivalent to some fibre bundle
$(E,\pi,A\times M)$
over the composite base $A\times M$. In fact, if
$\{(E_a,\pi_a,M),\ a\in A\}$
is given, we construct the fibre bundle $(E,\pi,A\times M)$ by putting
\[
E=\bigcup_{a\in A}E_a,\qquad
\pi:E\to A\times M,\quad \pi(u)=(a_u,\pi_{a_u}(u)),
\]
where  $u\in E$ and $a_u$ is the unique  $a_u\in A$ for which
$E_{a_u}\ni u$.
Conversely, if $(E,\pi,A\times M)$ is given, we can construct
$\{(E_a,\pi_a,M),\ a\in A\}$ through
\[
E_a=\bigcup_{x\in M}\pi^{-1}(a,x),\quad\pi_a:E_a\to M,\quad\pi_a(u)=x_u,
\]
where $a\in A,\ u\in E_a,\mbox{ and } x_u$ is the unique $x_u\in M$ for
which $\pi^{-1}(a,x_u)\ni u$.

	Now it is trivial to check that the above-defined transports are
 ${^{a,b}}I_{x\to y} = K_{(a,x)\to(b,y)}^{id_{A\times M}}$,
i.e. they are equivalent to the (flat) transport along the identity map
$id_{A\times M}:A\times M\to A\times M$ in  $(E,\pi,A\times M)$.
	\end{Exam}

	\begin{Exam}{Exam3.2}
	The described in the previous example general construction can be
specified in the case of a flat transport in the fibre bundle of
tensors of a fixed rank $k\in\mathbb{N}\cup\{0\}$ over a differentiable
manifold $M$ as follows.

	Let
$A_r = \{(p,q):\quad p,q\in\mathbb{N}\cup\{0\},\quad p+q=r\}$ and
$\left.T_{q}^{p}\right|_x(M)$,  $(p,q)\in A_r$ be the tensor space
of type $(p,q)$ over $x\in M$.
	The \textbf{tensor bundle of type $\mathbf{(p,q)}$}
is $\xi_{(p,q)} := \left(T_{q}^{p}(M),\pi_{(p,q)},M\right)$ with
 $T_{q}^{p}(M) := \bigcup_{x\in  M}\left.T_{q}^{p}\right|_x(M)$ and
 $\pi_{(p,q)}(u)=x, \mbox{ for } u\in \left.T_{q}^{p}\right|_x(M)$,
 $x$ being the unique $x\in M$ for which
$\left.T_{q}^{p}\right|_x(M)\ni u$.
  	The \textbf{tensor bundle
$(T_r(M),\pi_r,A_r\times M)$ of rank $\mathbf{r}$ }
is constructed by the above scheme:
 $T_r(M) := \bigcup_{(p,q)\in A_r}T_{q}^{p}(M) =
\bigcup_{p+q=r}\bigcup_{x\in M}\left.T_{q}^{p}\right|_x(M)$ and
\(
\pi_r(u) = \left((p,q),x\right)
\)
for  $u\in T_r(M)$ with  $(p,q)\in A_r$ and
$x\in M$ being defined by  $\left.T_{q}^{p}\right|_x(M)\ni u$.
So
\(
\pi_{r}^{-1}((p,q),x) =
\pi_{(p,q)}^{-1}(x) = \left.T_{q}^{p}\right|_x(M).
\)
Then the desired transport along $id_{A_r\times M}$ is described by
the maps
 ${^r\!}K_{((p',q'),x')\to ((p'',q''),x'')}^{id_{A\times M}}$,
 $p'+q'=p''+q''=r,\ x',x''\in M $
satisfying the relations:
\[
\begin{array}{c}
{^r\!}K_{((p',q'),x')\to((p'',q''),x'')}^{id_{A_r\times M}}:
\left.T_{q'}^{p'}\right|_{x'}(M) \to
				\left.T_{q''}^{p''}\right|_{x''}(M),
\\[1.8ex]
{^r\!}K_{((p'',q''),x'')\to((p''',q'''),x''')}^{id_{A_r\times M}} \circ
{^r\!}K_{((p',q'),x')\to((p'',q''),x'')}^{id_{A_r\times M}} \>=
\qquad\qquad\qquad\qquad \qquad\qquad  \\[1.3ex]
\qquad\qquad\qquad \qquad\qquad\qquad\qquad \qquad\qquad  =\>
{^r\!}K_{((p',q'),x')\to((p''',q'''),x''')}^{id_{A_r\times M}},
\\[1.8ex]
{^r\!}K_{((p',q'),x')\to((p',q'),x')}^{id_{A_r\times M}} =
id_{\left.T_{q'}^{p'}\right|_{x'}(M)},
\\[1ex]
\qquad p'+q'=p''+q''=p'''+q'''=r, \qquad x',x'',x'''\in M.
\end{array}
\]
	\end{Exam}

	\begin{Exam}{Exam3.3}
	Every transport along $\varkappa'':M\to M'$ in a bundle
$(E_0,\pi_0,M')$ induces a transport along
$\varkappa'\times\varkappa''$ for
$\varkappa':A\to A'$ in any bundle  $(E,\pi,A'\times M')$ for which the
fibres $\pi^{-1}(a',x')$ are homeomorphic to $\pi_{0}^{-1}(y')$ for any
$a'\in A'\mbox{ and }x',y'\in M'$. In fact, by example~\ref{Exam3.1},
$(E,\pi,A'\times M')$ is equivalent to the family
$\{\xi_{a'}=(E_{a'},\pi_{a'},M'),\quad a'\in A'\}$ with
$E_{a'}=\bigcup_{x'\in M'}\pi^{-1}(a',x')$ and  $\pi_{a'}(u)=x'$ for
$u\in\pi^{-1}(a',x')$. As the fibres of all the introduced fibre
bundles are homeomorphic there are fibre morphisms
$(h_{a'},id_{M'})$ from $\xi_{a'}$ on
$\xi_0$, i.e.  $h_{a'}:E_{a'}\to E_0$ and
$\pi_{a'}=\pi_0\circ h_{a'}$,
$a'\in A'$. Then it is easy to verify that the maps
	\begin{eqnarray*}
\lefteqn{
K_{(a,x)\to(b,y)}^{\varkappa'\times\varkappa''}  :=
\left(
\left.h_{\varkappa'(b)}\right|_{\pi^{-1}(\varkappa'(b),\varkappa''(y))}
\right)^{-1} \circ
} \\ & & \circ
\left(
\left.
h_{\varkappa'(a)}\right|_{\pi^{-1}(\varkappa'(a),\varkappa''(x))}
\right)
:\pi^{-1}(\varkappa'(a),\varkappa''(x)) \to
				\pi^{-1}(\varkappa'(b),\varkappa''(y))
	\end{eqnarray*}
define a transport along  $\varkappa'\times\varkappa''$ in
$(E,\pi,A'\times M')$.
	\end{Exam}

	\begin{Exam}{Exam3.4}
	This example is analogous to example~\ref{Exam3.2} and is
obtained from it by replacing  $p$ and  $q$
by integer functions over $M$.

	Let
$f,g:M\to\mathbb{N}\cup\{0\},\ f+g=r\in\mathbb{N}\cup\{0\}$, and
\[
\xi_{(f,g)}^{r}:=\left(^r\!\mathcal{T}_{g}^{f}(M),\pi_{f,g},M\right),
\quad ^r\!\mathcal{T}_{g}^{f}(M):=
\bigcup_{x\in M}\!\left.T_{g(x)}^{f(x)}\!\right|_x(M)
\mbox{ and } \pi_{f,g}(u):=x
\]
for $u\in\left.T_{g(x)}^{f(x)}\right|_x(M)$.
The transports along maps in the so defined fibre bundle preserve the
tensor's rank, but, generally, they change the tensor's type.
	\end{Exam}

\section {\bfseries Linear transports along maps}
\label{IV}
\setcounter {equation} {0}

	In this section $\xi:=(E,\pi,B)$ is supposed to be a complex (or
real) \emph{vector} bundle.

	As we said above, a
\emph{linear transport (or L-transport) along maps}
in a vector bundle is one satisfying eq.~(\ref{2.10}).
For these transports \emph{mutatis mutandis}
valid are almost all definitions and
results concerning linear transports along paths in (vector) fibre
bundles~\cite{f-LTP-general,f-LTP-appl,f-LTP-Cur+Tor,f-LTP-metrics}.
This is true for the cases in which the fact
that a path  is a map
from a real interval into some set is not explicitly used.
In particular, by replacing
the path $\gamma:J\to B$ and the linear transport along paths
$L_{s\to t}^{\gamma}$, $s,t\in J$ with, respectively, a map
$\varkappa:N\to B$
and a linear transport along maps
$L_{l\to m}^{\varkappa}$, $l,m\in N$, one
obtains a valid version of sections 2 and 3 of~\cite{f-LTP-general},
section 3 to proposition 3.3 and section 5 to proposition 5.3
of~\cite{f-LTP-appl}, and sections 1, 2, and 4 of~\cite{f-LTP-metrics}.
The other parts of these works, as well as~\cite{f-LTP-Cur+Tor},
deal more or less with explicit properties
of the real interval $J$, mainly via the differentiation along paths
$\mathcal{D}^\gamma$ (or $\mathcal{D}^\gamma_s$)~\cite{f-LTP-general}.
These exceptional definitions and results can, if possible, be
generalized as follows.

	Let $\mathbf{N}$ be a neighborhood in $\mathbb{R}^k$,
$k\in\mathbb{N}$, e.g. one may take
$\mathbf{N} = J\times\cdots\times J$
($k$-times), $J$ being a real interval.
So, any $l\in\mathbf{N}$ has a form
$l=(l^1,\ldots,l^k)\in\mathbb{R}^k$. We put
$\varepsilon_a:=(0,\ldots,0,\varepsilon,0,\ldots,0)\in\mathbb{R}^k$
where
$\varepsilon\in\mathbb{R}$ stands in the
$a$-th position, $1\le a\le k$.

	Let  $\mathrm{Sec}^p(\xi)$ (resp. $\mathrm{Sec}(\xi)$) be the
set of $C^p$ (resp. all) sections over $\xi$ and
$L_{l\to m}^{\varkappa}$ be a $C^1$
(on $l$) linear transport along $\varkappa:\mathbf{N}\to B$.
Now definition 4.1 from~\cite{f-LTP-general} is replaced by

	\begin{Def}	\label{Def4.1}
	The $a$-th, $1\le a\le k$ partial derivation along maps
generated by  $L$ is a  map
${_a  }\mathcal{D}:\varkappa\mapsto{_a  }\mathcal{D}^\varkappa$
where the $a$-th (partial) derivation along $\varkappa$
(generated by $L$) ${_a  }\mathcal{D}^\varkappa$ is  a map
	\begin{equation} \label{4.1}
{_a  }\mathcal{D}^\varkappa :
\mathrm{Sec}^1\left(\left.\xi\right|_{\varkappa(\mathbf{N})}\right) \to
\mathrm{Sec}\left(\left.\xi\right|_{\varkappa(\mathbf{N})}\right)
	\end{equation}
defined for
$\sigma\in\mathrm{Sec}^1\left(\xi|_{\varkappa(\mathbf{N})}\right)$ by
	\begin{equation} 	\label{4.2}
\left({_a  }\mathcal{D}^\varkappa\sigma\right) (\varkappa(l)) :=
\lim_{\varepsilon\to 0}  	\left[ \frac{1}{\varepsilon} \left(
L_{l+\varepsilon^a \to l}^{\varkappa}\sigma(\varkappa(l+\varepsilon^a))
				- \sigma(\varkappa(l)) \right) \right].
	\end{equation}
	The (partial)  $a$-th derivative of $\sigma$ along $\varkappa$
with respect to  $L$ is ${_a  }\mathcal{D}^\varkappa\sigma$.
Its value at $\varkappa(l)$ is
given by the operator
	\begin{equation}	\label{4.3}
	\mathcal{D}_{l^a}^{\varkappa} :
\mathrm{Sec}^1\left(\left.\xi\right|_{\varkappa(\mathbf{N})}\right) \to
\pi^{-1}(\varkappa(l))
	\end{equation}
by
$\mathcal{D}_{l^a}^{\varkappa}\sigma :=
\left({_a  }\mathcal{D}^{\varkappa}\sigma\right) (\varkappa(l))$.
	\end{Def}

	Evidently, for $k=1$ this definition reduces to definition 4.1
of~\cite{f-LTP-general}.

	On the basis of the above definition, almost all of the
above-mentioned exceptional definitions and results can be modified
by replacing in them
	 $\gamma:J\to B$,  $L_{s\to t}^{\gamma}$, $\mathcal{D}^\gamma$,
 and $\mathcal{D}_{s}^{\gamma}$, respectively with
 $\varkappa:\mathbf{N}\to B$,
$L_{l\to m}^{\varkappa}$, ${_a  }\mathcal{D}^\varkappa$,
and ${_a  }\mathcal{D}_{l^a}^{\varkappa}$ .
Below we sketch some results in this field.

	A corollary of~(\ref{2.2}),~(\ref{2.3}), and~(\ref{4.2}) is

	\begin{Prop}	\label{Prop4.1}
	The operators ${_a  }\mathcal{D}^\varkappa$ are
($\mathbb{C}$-)linear
and
	\begin{equation}	\label{4.4}
\mathcal{D}_{m^a}^{\varkappa} \circ L_{l\to m}^{\varkappa} \equiv 0 .
	\end{equation} \end{Prop}

	If $\{e_i\}$ is a field of bases on $\varkappa(\mathbf{N})$,
i.e. $\{e_i(n)\}$ is a basis in $\pi^{-1}(\varkappa(n))$,
$n\in\mathbf{N}$,
then in it the L-transport $L$ along maps is described by the matrix
\(
H(m,l;\varkappa)= \left[H_{\;i}^{j}(m,l;\varkappa)\right]_{i,j=1}^{n},
\)
$n:=\dim_\mathbb{C}(\pi^{-1}(\varkappa(l))$, which is defined by
\(
L_{l\to m}^{\varkappa}e_i(l) =: H_{\;i}^{j}(m,l;\varkappa) e_j(m).
\)
	A simple calculation of the limit in~(\ref{4.2}) verifies

	\begin{Prop}	\label{Prop4.2}
	If
$\sigma\in\mathrm{Sec}^1\left(\xi|_{\varkappa(\mathbf{N})}\right)$,
then
	\begin{equation}	\label{4.5}
	\mathcal{D}_{l^a}^{\varkappa} \sigma =
\left[
\frac{\partial\sigma^i(\varkappa(l))}{\partial l^a} +
{_a  }\Gamma_{\;j}^{i}(l;\varkappa)\sigma^j(\varkappa(l))
\right] e_i(l) ,
	\end{equation}
where the \textbf{components} of  $L$ are defined by
	\begin{equation}	\label{4.6}
{_a  }\Gamma_{\;j}^{i}(l;\varkappa) :=
\left.
\frac{\partial H_{\;j}^{i}(l,m;\varkappa)}{\partial m^a}
\right|_{m=l} =
- \left.
\frac{\partial H_{\;j}^{i}(m,l;\varkappa)}{\partial m^a}
\right|_{m=l}.
	\end{equation}
	\end{Prop}

	The components of $L$ satisfy
	\begin{equation}	\label{4.7}
{_a  }\mathcal{D}^\varkappa e_j =
			\left({_a  }\Gamma_\varkappa\right)_{\;j}^{i}e_i
	\end{equation}
and form $k$ matrices
${_a  }\Gamma_\varkappa(l) :=
\left[{_a  }\Gamma^{i}_{\;j}(l;\varkappa)\right]_{i,j=1}^{n}$,
 $a=1,\ldots,k$ which under the transformation
 $e_j(l)\mapsto e_{j}^{\prime}(l) = A_{j}^{i}(l)e_i(l)$    change to
	\begin{equation}	\label{4.8}
{_a  }\Gamma_\varkappa^\prime(l) =
A^{-1}(l) \left({_a  }\Gamma_\varkappa(l)\right) A(l) +
A^{-1}(l)\frac{\partial}{\partial l^a} A(l)
	\end{equation}
with $A(l):=\left[A_{j}^{i}(l)\right]$, which is a simple corollary
of~(\ref{4.7}).
Hence, the difference of the matrices  ${_a  }\Gamma_\varkappa$
of two L-transports along one and the same map behaves like a tensor
of type $(1,1)$ under a transformation of the bases.

	On the above background one can \emph{mutatis mutandis}
reformulate the remaining part of Sect.~4 of~\cite{f-LTP-general}.
In particular, in this way
is established the equivalence of the sets of L-transports along maps
$\varkappa:\mathbf{N}\to B$, $\mathbf{N}\subseteq\mathbb{R}^k$ and the
one of partial derivations along maps. Sect.~6 and
the rest of Sect.~3 and Sect.~5 of~\cite{f-LTP-appl}
can be modified analogously, only in the last case the tangent vector
field $\dot\gamma$ to $\gamma:J\to M$ has to be replaced with
the set of tangent vectors $\{\dot\varkappa_a\}$ to $\varkappa$,
\(
\dot\varkappa_a(l) :=
\left(\frac{\partial\varkappa^i(l)}{\partial l^a}\right)
\left(
\left. \frac{\partial} {\partial x^i}\right|_{\varkappa(l)}
\right).
\)

	The introduction of torsion and curvature needs more details
which will be presented below.

	Let $M$ be a differentiable manifold and there be given a
$C^1$ map $\eta:N\times N'\to M$, with $N$ being a neighborhood in
$\mathbb{R}^k$
and  $N'$ - in $\mathbb{R}^{k'}$, $k,k'\in\mathbb{N}$. Let
$\eta(\cdot,m):l\mapsto\eta(l,m)$, $\eta(l,\cdot):m\mapsto\eta(l,m)$
for  $(l,m)\in N\times N^\prime$.
Let
$\eta_{a}^{\prime}(\cdot,m)$, and  $\eta_{b}^{\prime\prime}(l,\cdot)$,
 $l\in N$, $m\in N'$, $a=1,\ldots,k$, $b=1,\ldots,k'$ be the  tangent
vector fields to  $\eta(\cdot,m)$ and $\eta(l,\cdot)$, respectively.

	\begin{Def}	\label{Def4.2}
	The torsion operators of an L-transport along maps in the
tangent bundle  $(T(M),\pi,M)$ are maps
\(
\mathcal{T}_{a,b}:\eta\mapsto\mathcal{T}_{a,b}^{\eta}:N\times N'\to T(M)
\)
which for $(l,m)\in N\times N'$ are given by
	\begin{equation}	\label{4.9}
\mathcal{T}_{a,b}^{\eta}(l,m) :=
\mathcal{D}_{l^a}^{\eta(\cdot,m)} \eta_{b}^{\prime\prime}(\cdot,m) -
\mathcal{D}_{m^b}^{\eta(l,\cdot)} \eta_{a}^{\prime}(l,\cdot)
\in T_{\eta(l,m)}(M) .
	\end{equation}
	\end{Def}

	Similarly, for $\eta:N\times N'\to B$, $B$ being the base of
 a vector bundle $(E,\pi,B)$, we have

	\begin{Def}	\label{Def4.3}
	The curvature operators of an L-transport along maps in
$(E,\pi,B)$ are maps
\[
\mathcal{R}_{a,b}:\eta\mapsto\mathcal{R}_{a,b}^{\eta}:(l,m)\mapsto
\mathcal{R}_{a,b}^{\eta}(l,m):\mathrm{Sec}^2(E,\pi,B) \to
\mathrm{Sec}(E,\pi,B)
\]
defined for $(l,m)\in N\times N'$ by
	\begin{equation}	\label{4.10}
\mathcal{R}_{a,b}^{\eta}(l,m) :=
{_a  }\mathcal{D}^{\eta(\cdot,m)} \circ
			{_b\!}\mathcal{D}^{\eta(l,\cdot)} -
{_b\!}\mathcal{D}^{\eta(l,\cdot)} \circ
			{_a  }\mathcal{D}^{\eta(\cdot,m)}.
	\end{equation}
	\end{Def}

	The further treatment of curvature and torsion can be done
by the same methods as in~\cite{f-LTP-Cur+Tor,f-LTP-Cur+Tor-prop}
(cf.~\cite[Sect.~8]{f-Jacobi}).

	In the composite case there arises a kind of `restricted'
partial derivation along maps generated by L-transports along maps.

	Let $N=A\times M$ with $M$ being a neighborhood in
$\mathbb{R}^k$ for some $k\in\mathbb{N}$.

	In this case instead of definition~\ref{Def4.1}, we have

	\begin{Def}	\label{Def4.4}
	The $a$-th, $1\le a\le k$, partial derivative of type $\beta,\
\beta\in A$ along the map  $\varkappa:A\times M\to B$ generated by an
L-transport $L$ along maps in a vector bundle $\xi=(E,\pi,B)$ is a map
\(
{_{a}^{\beta}\!}\mathcal{D}:\varkappa \mapsto
{_{a}^{\beta}\!}\mathcal{D}^\varkappa
\),
where
	\begin{equation}	\label{4.11}
{_{a}^{\beta}\!}\mathcal{D}^\varkappa
:\mathrm{Sec}^1\left( \left.\xi\right| _{\varkappa(A\times M)} \right)
\to \mathrm{Sec}\left( \left.\xi\right| _{\varkappa(A\times M)} \right)
	\end{equation}
is the partial derivation along  $\varkappa$ (generated by $L$)
which is defined by the equation
	\begin{eqnarray}	\nonumber
\left(
{_{a}^{\beta}\!}\mathcal{D}^\varkappa \sigma
\right)
(\varkappa(\alpha,x))
 & := &
\lim_{\varepsilon\to 0} \Biggl[ \frac{1}{\varepsilon} \Bigl(
L_{(\alpha,x+\varepsilon^a) \to (\beta,x)}^{\varkappa}
\sigma(\varkappa(\alpha,x+\varepsilon^a)) -
\\   	\label{4.12}	& -  &
L_{(\alpha,x)\to(\beta,x)}^{\varkappa}   \sigma(\varkappa(\alpha,x))
\Bigr) \Biggr]
	\end{eqnarray}
for
\(
\sigma\in\mathrm{Sec}^1\left(\left.\xi\right|_{\varkappa(A\times M)}
\right)
\)
and  $(\alpha,x)\in A\times M$.
The $a$-th (partial) derivative (of type $\beta$) of  $\sigma$ along
$\varkappa$ with respect of $L$ is
${_{a}^{\beta}\!}\mathcal{D}^\varkappa \sigma$.
Its value at $\varkappa(\alpha,x)$ is given by the operator
	\begin{equation}	\label{4.13}
{^{\beta}\!}\mathcal{D}_{(\alpha,x^a)}^\varkappa :
\mathrm{Sec}^1\left(\left.\xi\right|_{\varkappa(A\times M)} \right) \to
\pi^{-1}(\varkappa(\beta,x))
	\end{equation}
by
\(
{^{\beta}\!}\mathcal{D}_{(\alpha,x^a)}^\varkappa \sigma :=
\left(
{_{a}^{\beta}\!}\mathcal{D}^\varkappa \sigma
\right)
\varkappa(\alpha,x),
\)
with $x^a$ being the $a$-th component of $x\in M\subseteq\mathbb{R}^k$.
	\end{Def}

	For $A=\emptyset$ this definition reduces to
definition~\ref{Def4.1}.

	Notice that the operator
 ${_{x\!}^{a,b}}\nabla_{V}^{I}$ used in~\cite[eq.~(7.14)]{f-Jacobi}
is a special case of
  ${_{a}^{\beta}\!}\mathcal{D}^\varkappa$, viz.
\(
\left( {_{x}^{a,b}\!}\nabla_{V}^{I} (\sigma) \right) (a,x) =
\sum_{i} V^i \left({_{i}^{b}\!}\mathcal{D}^{id_M} \sigma \right) (a,x)
\)
for $i=1,\ldots,\dim M$, $ a,b\in A,\ x\in M$, $M$ being a differentiable
manifold and $V$ being a vector field on $M$ with local components $V^i$.

	Now the corresponding results
from~\cite{f-LTP-general,f-LTP-appl,f-LTP-Cur+Tor,f-LTP-metrics} can be
modified step by step on the basis of definition~\ref{Def4.4} in the
above-described way, where definition~\ref{Def4.1} was used.

\section {\bfseries Consistency with bundle morphisms}
\label{V}
\setcounter{equation}{0}

	The work~\cite{f-TP-morphisms} investigates problems concerning
the consistency of transports along paths in fibre bundles and bundle
morphisms between them. A critical reading of this paper reveals the
insignificance of the fact that the transports in it are along paths;
nowhere there is the fact used that the path $\gamma$ is a map from a
real interval  $J$ into the base $B$ of some fibre bundle. For this
reason all of the work~\cite{f-TP-morphisms} is valid
\emph{mutatis mutandis} for arbitrary transports
along maps; one has simply to replace the transports along paths, like
$I_{s\to t}^{\gamma}$, $\gamma:J\to B$, $s,t\in J$, with transports
along arbitrary maps, like $K_{l\to m}^{\varkappa}$,
$\varkappa:N\to B$, $l,m\in N$.
Below are stated \emph{mutatis mutandis} only some definitions and
results from~\cite{f-TP-morphisms}. There proofs are omitted as they
can  easily be obtained from the corresponding ones
in~\cite{f-TP-morphisms}.

	Let there be given two fibre bundles
$\xi _{h}:=(E_{h},\pi _{h},B_{h}),\quad h=1,2$ in
which defined are, respectively, the transports along maps ${}^1 \! K$
and
${}^2 \! K$. Let $(F,f)$ be a bundle morphism from
$\xi _{1}$ into $\xi _{2}$,
i.e.
$F:E_{1} \to E_{2}, \quad f:B_{1} \to B_{2}$ and
$\pi _{2}\circ F=f\circ \pi _{1}$~\cite{Husemoller} .
Let $F_{x}:=\left.F\right|_{\pi ^{-1}_{1}(x)}$ for $x\in B_{1}$ and
$\varkappa:N \to B_{1}$ be an arbitrary map in $B_{1}$.

	\begin{Def} \label{Def5.1}
	The bundle morphism $(F,f)$ and the pair $(^{1\!}K,{^2\!}K)$ of
transports, or the transports ${^1\!}K$ and ${^2\!}K$, along maps are
consistent (resp. along the map $\varkappa$) if they commute in a sense
that the equality
	\begin{equation}	\label{5.1}
F_{\varkappa(m)}\circ {^1\!}K^{\varkappa}_{l\to m}=
{^2\!}K^{f\circ\varkappa}_{l \to m}\circ F_{\varkappa(l)},
\quad l,m\in N
\end{equation}
is fulfilled for every (resp. the given) map $\varkappa$.
	\end{Def}

	A special case of definition~\ref{Def5.1} is the
condition~(\ref{2.9}) for consistency with a bundle binary operation
(in particular, a bundle metric), which is obtained from it for:
 $\xi_1=(E,\pi,B)\times(E,\pi,B)$,  $\xi_2=(M,\pi_0,m_0)$ with a fixed
$m_0\in M$, $\pi_0:M\to\{m_0\}$, $F_x=\beta_x$ with $x\in B$,
$f:B\times B\to\{m_0\}$,
\(
{^1\!}K^{\varkappa}_{l\to m} =
K^{\varkappa}_{l\to m} \times K^{\varkappa}_{l\to m},
\)
and  ${^2\!}K^{f\circ\varkappa}_{l \to m} = id_M$.

	\par
	Let, in accordance with theorem~\ref{Th2.1}
	(cf.~\cite[theorem 3.1]{f-TP-general}),
there be chosen sets $Q_{1}$ and $Q_{2}$ and one-to-one maps
\(
{^h\!}F^{\varkappa_{h}}_{l_{h}}:\pi ^{-1}_{h}(\varkappa_{h}(l_{h}))
\to Q_{h},
\)
$h=1,2$, which are associated, respectively, with the maps
\(
\varkappa_{h}:N_{h} \to B_{h},
\quad l_{h}\in N_{h}, \quad h=1,2
\)
and are such
that (cf.~(\ref{2.4}))
	\begin{equation}	\label{5.2}
{^h\!}K^{\varkappa_{h}}_{l_{h}\to m_{h}}=
\left( {^h\!}F^{\varkappa_{h}}_{m_{h}}  \right) ^{-1} \circ
{^h\!}F^{\varkappa_{h}}_{l_{h}}, \qquad l_{h}, \ m_{h}\in N_{h},
\quad h=1,2.
\end{equation}

	\begin{Prop} \label{Prop5.1}
	The bundle morphism $(F,f)$ and the pair
$({^1\!}K,{^2\!}K)$ of tran\-s\-ports along maps, which are given
 by~(\ref{5.2}) by means of the maps ${^1\!}F$ and ${^2\!}F$, are
 consistent (resp. along a map $\varkappa$) iff there exists a map
	\begin{equation} 	\label{5.3}
 C_{0}(\varkappa,f\circ\varkappa):Q_{1} \to Q_{2},
	\end{equation}
 such that
	\begin{equation}	\label{5.4}
F_{\varkappa(l)} =
		\left( {^2\!}F^{f\circ\varkappa}_{l} \right) ^{-1}\circ
C_{0}(\varkappa,f\circ\varkappa)\circ
		\left( {^1\!}F^{\varkappa}_{l} \right),
	\end{equation}
or, equivalently, that
	\begin{equation}	\label{5.5}
F_{\varkappa(l)}={^2\!}K^{f\circ\varkappa}_{l_0\to m} \circ
C(l_{0};\varkappa,f\circ\varkappa) \circ {^1\!}K^{\varkappa}_{l\to l_0},
	\end{equation}
where $l_{0}\in N$ is arbitrary and
	\begin{equation}	\label{5.6}
C(l_{0};\varkappa,f\circ\varkappa) :=
\left( {^2\!}F^{f\circ\varkappa}_{l_0} \right) ^{-1} \circ
C_{0}(\varkappa,f\circ\varkappa) \circ
\left( {^1\!}F^{\varkappa}_{l_0} \right)
	\end{equation}
for every (resp. the given) map $\varkappa$.
	\end{Prop}

	\par
	Let there be given two fibre bundles
$\xi _{h}=(E_{h},\pi _{h},B_{h})$, $h=1,2$. We define the
\emph{
fibre bundle
$\xi _{0}=(E_{0},\pi _{0},B_{1})$ of bundle morphisms
}
from $\xi _{1}$ onto $\xi _{2}$ in the following way:
	\begin{eqnarray}
&\!\! \!\! \!\! \!\! \!\! \!
E_{0}:=\{ (F_{b_{1}},f)\! :\> F_{b_{1}}\! \! :\pi ^{-1}_{1}(b_{1}) \to
	\pi ^{-1}_{2}(f(b_{1})), \  b_{1}\in B_{1}, \
		f:B_{1} \to B_{2} \}, \quad   & \label{5.7}\\
 & \pi _{0}((F_{b_{1}},f)):=b_{1}, \quad (F_{b_{1}},f)\in E_{0}, \quad
 b_{1}\in B_{1}. &	\label{5.8}
	\end{eqnarray}

	It is clear that every section
$(F,f)\in \mathrm{Sec}(\xi _{0})$
is a bundle morphism from $\xi _{1}$ into $\xi _{2}$ and vice versa,
every bundle morphism from $\xi _{1}$
onto $\xi _{2}$ is a section of $\xi _{0}$.
(Thus a bundle structure in the
set $\mathrm{Morf}(\xi _{1},\xi _{2})$ of bundle morphisms from
 $\xi _{1}$ on $\xi _{2}$ is introduced.)

	If in $\xi _{0}$ a transport $K$ along the maps
in $B_{1}$ is given, then, according to definition~\ref{Def2.2}
(see eq.~(\ref{2.6})), the bundle morphism
$(F,f)\in \mathrm{Sec}(\xi _{0})$ is
(\mbox{$K$-})transported along $\varkappa:N\to B_{1}$ if
	\begin{equation}	\label{5.9}
 (F_{\varkappa(m)},f) =
	K^{\varkappa}_{l\to m}(F_{\varkappa(l)},f), \qquad l,m\in N.
	\end{equation}

	Given in  $\xi _{1}$ and $\xi _{2}$ the respective
transports ${^1\!}K$ and ${^2\!}K$ along the maps in $B_{1}$ and $B_{2}$
respectively. They generate
in $\xi _{0} $ a `natural' transport ${^0\!}K$ along the maps in $B_{1}$.
The action of this transport along
$\varkappa:N\to B_{1}$ on
$(F_{\varkappa(l)},f)\in \pi ^{-1}_{0}(\varkappa(l))$
for a fixed $l\in N$ and arbitrary $m\in N$ is defined by
	\begin{equation}	\label{5.10}
 {^0\!}K^{\varkappa}_{l\to m}(F_{\varkappa(l)},f) :=
 	\left( {^2\!}K^{f\circ\varkappa}_{l\to m} \circ F_{\varkappa(l)}
\circ {^1\!}K^{\varkappa}_{m\to l}, \; f \right) \in
					\pi ^{-1}_{0}(\varkappa(m)).
	\end{equation}
	\begin{Lem}             \label{Lem5.1}
If $(F,f)\in \mathrm{Sec}(\xi _{0})$, then~(\ref{5.1}) is  equivalent to
	\begin{equation}      	\label{5.11}
(F_{\varkappa(m)},f)={^0\!}K^{\varkappa}_{l\to m}(F_{\varkappa(l)},f),
	\qquad l,m\in M.
	\end{equation}
	\end{Lem}

	\begin{Prop}	\label{Prop5.2}
	The bundle morphism $(F,f)$ and the pair $({^1\!}K,{^2\!}K)$ of
tran\-s\-ports along maps are consistent (resp. along the map
$\varkappa$) if and only if $(F,f)$ is transported along every
(resp. the given) map $\varkappa$ with the help of the defined
from $({^1\!}K,{^2\!}K)$ in $\xi _{0}$ transport along maps ${^0\!}K$.
	\end{Prop}

\section {\bfseries Concluding discussion}
\label{VI}
\setcounter {equation} {0}

	\vspace{4ex}\noindent\hspace{4em}%
{\bfseries\large\sffamily\upshape(1)}
	The substitution of~(\ref{3.8}) into~(\ref{3.7}) gives
	\begin{equation}	\label{6.1}
K_{(a,x)\to(b,y)}^{\eta} =
\left( {^b\!}G_{y}^{\eta} \right)^{-1} \circ C_{a\to b}^{\eta} \circ
\left( {^a\!}G_{x}^{\eta} \right) =
\left( {^y\!}H_{b}^{\eta} \right)^{-1} \circ D_{x\to y}^{\eta} \circ
\left( {^x\!}H_{a}^{\eta} \right)
	\end{equation}
where (cf. (\ref{3.4}))
	\begin{equation}
C_{a\to b}^{\eta} :=
\left( C_{b}^{\eta} \right)^{-1} \circ C_{a}^{\eta} :Q_G\to Q_G,\quad
D_{x\to y}^{\eta} :=
\left( D_{y}^{\eta} \right)^{-1} \circ D_{x}^{\eta} :Q_H\to Q_H
	\end{equation}
are `transport like' maps. From them transports along
the identity map in corresponding fibre bundles can be constructed.
For instance, for
 $D_{x\to y}^{\eta}$ this can be done as follows.
	Consider the fibre bundle $(M\times Q_H,\pi_1,M)$ with
$\pi_1(x,q):=x$, $(x,q)\in M\times Q_H$. Hence
$\pi_{1}^{-1}(x)=\{x\}\times Q_H$. Defining
 $P_x : \{x\}\times Q_H\to Q_H$ by $P_x(x,q):=q$, we see that
 \(
\overline{D}_{x\to y}^{\,\eta} :=
P_{y}^{-1} \circ D_{x\to y}^{\eta} \circ P_x
\)
is a transport along $id_M$ in $(M\times Q_H,\pi_1,M)$. It depends on
$\eta$ as on a parameter. Consequently, we can write
	\begin{equation}	\label{6.3}
K_{(a,x)\to(b,y)}^{\eta} =
\left( P_{y}^{-1}\circ{^y\!}H_{b}^{\eta} \right)^{-1}
	\circ 	\overline{D}_{x\to y}^{\,\eta} 	\circ
\left( P_{x}^{-1}\circ{^x\!}H_{a}^{\eta} \right) .
	\end{equation}

	This decomposition is important when
$P_{x}^{-1}\circ{^x\!}H_{a}^{\eta}$
is independent of $x\in M$. Such a situation is realized in the
special case when
$B=A'\times M'$ and  $\eta=\eta'\times\eta''$ with $\eta':A\to A'$ and
 $\eta'':M\to M'$, i.e. for (some) transports along
$\eta'\times\eta''$ in
$(E,\pi,A'\times M')$. According to example~\ref{Exam3.1} the
fibre bundle $(E,\pi,A'\times M')$ is equivalent to the family
\(
\{
\xi_{a'}:\quad\xi_{a'}=(E_{a'},\pi_{a'},M'), \qquad
E_{a'}=\bigcup_{x'\in M'}\pi^{-1}(a',x'), \quad
\pi_{a'}(u)=x', \mbox{ for } u\in\pi^{-1}(a',x'), \quad a'\in A'
\}
\).
	Let $\xi_0=(E_0,\pi_0,M')$ be any fibre bundle for which there
exist bundle morphisms $(h_{a'},id_{M'})$ from $\xi_{a'}$ into $\xi_0$,
i.e.
 $h_{a'}:E_{a'}\to E_0$ and  $\pi_{a'}=\pi_0\circ h_{a'},\ a\in A'$.
(The existence of $\xi_0$  and  $h_{a'}$ is a consequence from the
fact that the fibres of all the defined bundles are homeomorphic;
e.g. one may put
 $\xi_0=\xi_{b'}$ for some fixed $b'\in A'$.)

	If ${^0\!}K_{x\to y}^{\varkappa''}$ is any transport along
$\varkappa''$ in $\xi_0$, then a simple calculation shows that
	\begin{equation} 	\label{6.4}
K_{(a,x)\to(b,y)}^{\varkappa'\times\varkappa''} :=
h_{\varkappa'(b)}^{-1} \circ {^0\!}K_{x\to y}^{\varkappa''} \circ
h_{\varkappa'(a)}
	 \end{equation}
defines a transport along $\varkappa'\times\varkappa''$ in
$(E,\pi,A'\times M')$.
The opposite statement is, generally, not valid, i.e. not for every
transport along $\varkappa'\times\varkappa''$ in  the fibre bundle
$(E,\pi,A'\times M')$ there exists a decomposition like~(\ref{6.4}).

	\vspace{4ex}\noindent\hspace{4em}%
{\bfseries\large\sffamily\upshape(2)}
	In vector bundles, such as the tensor bundles over a
differentiable manifold, sometimes the problem arises of
comparing  or performing some operations with vectors from
different fibres, or speaking more freely, with vectors (defined)
at different points. A way for approaching such problems is
the following one.

	Let the fibre bundle $(E,\pi,B)$ be endowed with (maybe linear)
transport $K^0$ along the identity map of $B$ and, e.g., a binary
operation $\beta$,
$\beta:x\mapsto\beta_x:\pi^{-1}(x)\times\pi^{-1}(x)\to Q_x$
for some sets $Q_x$, $x\in B$.
The problem is to extend the operation $\beta$
on sets like $\pi^{-1}(y)\times\pi^{-1}(z)$,  $y,z\in B$. A possible
solution is to replace $\beta$ with
 $\{ \overline{\beta}_x \}$ for some maps:
\[
\overline{\beta}_x:(y,z)\mapsto
\overline{\beta}_{x}^{\, y,z}:\pi^{-1}(y)\times\pi^{-1}(z): \to Q_x,
\]
where
	\begin{equation}	\label{6.5}
	\overline{\beta}_{x}^{\, y,z} :=
\beta_x \circ \left( K_{y\to x}^{0} \times K_{z\to x}^{0} \right).
	\end{equation}

	For instance, if $Q_x=\pi^{-1}(x)$ and $(E,\pi,B)$ is a
vector bundle one can define in this way the linear combination of
vectors from different fibres by the equality
\[
(\lambda u + \mu v)_x :=
\lambda K_{y\to x}^{0}u + \mu K_{z\to x}^{0}v,\qquad
\lambda,\mu\in\mathbb{C},\ u\in\pi^{-1}(y),\ v\in\pi^{-1}(z).
\]
It depends on $x$ as on a parameter. If $K^0$ is linear, then
\(
(\lambda u + \mu v)_{x''} = K_{x'\to x''}^{0}(\lambda u + \mu v)_{x'}
\),
 $x',x''\in B$.

	Also in this way can be introduced different kinds of
integrations of the sections of  $(E,\pi,B)$ if on $Q_x$
corresponding measures are defined.
Viz. if $d\mu_x$ is a measure on $Q_x$ and
$\sigma\in\mathrm{Sec}(E,\pi,B)$, the integral of $\sigma$ over a set
$V\subseteq E$ is defined as
 $\int\limits_{V}^{}\left(K_{y\to x}^{0}\sigma(y)\right) d\mu_x$.
This  procedure is especially useful in tensor bundles in which there
are different possibilities depending on the understanding of
the product of the integrand with the measure, e.g. it can be
a tensor product that may be combined with some contraction(s) too.

	 The situation is important when $(E,\pi,B)$ is endowed with a
transport along maps of a given kind, i.e. along
$\varkappa\in\mathcal{K}$, where $\mathcal{K}$ is a certain set of
maps onto  $B$. A typical example of such a set is the set of all
paths on $B$, i.e.
 $\{\gamma:\quad \gamma:J\to B,\ J\subseteq\mathbb{R}\}$.

	Let for some $x\in B$ there be a neighborhood $U\ni x$ in
$B$ with the property that for any $y\in U$ there are a unique map
$\varkappa_y:N_y\to B$,  $\varkappa_y\in\mathcal{K}$ and a set
$M_y\subseteq N_y$
such that $\varkappa_y\left.\right|_{M_y}:M_y\to U$,
 $\varkappa_y\left.\right|_{M_y}(m_x)=x$, and
 $\varkappa_y\left.\right|_{M_y}(m_y)=y$ for some $m_x,m_y\in M_y$.
A well-known example of this kind is the case of geodesic paths (curves)
on a differentiable manifold endowed with an affine
connection~\cite{K&N-1,Postnikov-Morse}.

	In such a neighborhood $U$ one can repeat the above discussion
(of the flat case) with the only change that $K_{y\to x}^{0}$
has to be replaced with $K_{m_y\to m_x}^{\varkappa_y}$.

	The use of transports along the maps $\varkappa_y$ has the
disadvantage that the result depends on $\varkappa_y$, but as they
are unique in the above sense this is insignificant. If the fibre
bundle admits some `natural' family of such maps, as the above-pointed
case of geodesic curves, the question of this uniqueness does not arise
at all. If the set of maps with the considered property does not exists
or is not unique, then the pointed procedure does not exists or is not
unique and, consequently, one gets nothing or not a `reasonable' result,
respectively.

	\vspace{4ex}\noindent\hspace{4em}%
{\bfseries\large\sffamily\upshape(3)}
	The class of Hermitian (resp. real) metrics on a
complex (resp. real) differentiable manifold $M$ turns out
to be in one-to-one correspondence with the class of flat linear
transports in the tensor bundle of rank 1 over it
(see example~\ref{Exam3.2}). Below is presented the proof
of this statement.

	The tensor bundle of rank 1 over $M$ is
\[
\left(
T_{0}^{1}(M)\bigcup T_{1}^{0}(M),\pi_1,\{(1,0),(0,1)\}\times M
\right)
\]
with $\pi_1(u)=((p,q),x)$ for $u\in \left.T_{q}^{p}\right|_x(M)$,
$p+q=1$.
According to proposition~\ref{Prop3.1} there are the following four
kinds of transports along  $id_{ \{(1,0),(0,1)\}\times M }$:
	\begin{eqnarray}
L_{x\to y}^{} &:=&
L_{((1,0),x)\to((1,0),y)}^{id_{ \{(1,0),(0,1)\}\times M }}=
F_{y}^{-1}\circ F_x = G_{y}^{-1}\circ G_x, 	\label{6.6} \\
L_{x\to y}^{*} &:=&
L_{((0,1),x)\to((0,1),y)}^{id_{ \{(1,0),(0,1)\}\times M }}=
{^*\!}F_{y}^{-1}\circ {^*\!}F_x =
{^*\!}G_{y}^{-1}\circ {^*\!}G_x,\label{6.7} \\
L_{x\to y}^{1,0} &:=&
L_{((1,0),x)\to((0,1),y)}^{id_{ \{(1,0),(0,1)\}\times M}}
= {^*\!}F_{y}^{-1}\circ F_x =
{^*\!}G_{y}^{-1}\circ {^*\!}C^{-1}\circ C \circ G_x,	\label{6.8} \\
L_{x\to y}^{0,1} &:=&
L_{((0,1),x)\to((1,0),y)}^{id_{\{(1,0),(0,1)\}\times M }}
= F_{y}^{-1}\circ {^*\!}F_x =
G_{y}^{-1}\circ C\circ {^*\!}C \circ {^*\!}G_x.	\label{6.9}
	\end{eqnarray}
Here, for brevity, we have put:
\(
F_x:=F_{((1,0),x)}^{\ldots}=C\circ G_x,\
{^*\!}F_x:=F_{((0,1),x)}^{\ldots}={^*\!}C\circ {^*\!}G_x,\
G_x:={^{(1,0)} }G_{x}^{\ldots},\
{^*\!}G_x:={^{(0,1)} }G_{x}^{\ldots},\
C:=C_{(1,0)}^{\ldots},\
{^*\!}C:=C_{(0,1)}^{\ldots},\
\)
where the dots (\ldots)\ stand for ${id_{ \{(1,0),(0,1)\}\times M }}$.
These maps act between vector spaces and are linear because of the
linearity of the considered transports. Let in all vector spaces,
as well as in the fibres
$ T_x(M) := \left.T_{0}^{1}\right|_x(M)$ and
$ T_x^*(M) := \left.T_{1}^{0}\right|_x(M)$, $x\in M$,
some bases be fixed in which the matrix of a map
$X\in\{ F_x,\ ^*\!F_x,\ G_x,\ ^*\!G_x,\ C,\ ^*\!C,\ L_{x\to x}^{1,0}\}$
will be written as $[X]$.

	A (fibre) Hermitian (resp. real) metric on $M$ is
$g:x\mapsto g_x$,
$x\in M$. Here $g_x:T_x(M)\times T_x(M)\to\mathbb{C}$
(resp. $\mathbb{R}$) are
 $1\frac{1}{2}$ linear (resp. bilinear), nondegenerate, and Hermitian
(resp. symmetric) maps~\cite{Greub&et.al.-1}. Let in the above bases the
matrix of $g_x$ be  $G(x)$; we have $\det G(x)\not= 0$ and
$G^\dag(x)=G(x)$
(resp. $G^\top(x)=G(x)$ as in the real case $G^\dag=G^\top$),
where \dag\ (resp. $\top$) means Hermitian conjugation (resp.
transposition), i.e. transposition plus complex conjugation (denoted
by a bar). Because of this there is a unitary (resp. orthogonal) matrix
$D(x)$, i.e.
 $D^\dag=D^{-1}$ (resp. $D^\top=D^{-1}$), such that
 $G(x)=D^\dag(x)G_{p,q}D(x)$ with
\(
G_{p,q}:=\mathrm{diag}(\underbrace{+1,\ldots,+1}_{\mathrm{p-times}},
			\underbrace{-1,\ldots,-1}_{\mathrm{q-times}})
\)
for some unique $p,q\in\mathbb{N}\bigcup\{0\}$,
$p+q=\dim M$~\cite{Bellman}.

	Now the idea is to interpret the maps~(\ref{6.8}) for $y=x$ as
metrics.

	In fact, if for $u,v\in T_x(M)$ we put
 $g_x(u,v):=(g_x(u,\cdot))\overline{v}$ with
	\begin{equation}	\label{6.10}
	g_x(u,\cdot):=L_{x\to x}^{1,0}u \in T_{x}^{*}(M),
	\end{equation}
we find $G(x)=[{^*\!}F_x]^{-1}[F_x]$.
This matrix will be Hermitian if, e.g.,we choose
$[{^*\!}F_x]^{-1}=[F_x]^\dag$, i.e. ${^*\!}F_x^{-1}=F_x^\dag$,
which leads to
 $\left(L_{x\to x}\right)^\dag = L_{x\to x}^*$.
In particular, we can choose
${^*\!}C^{-1} = C^\dag$ and ${^*\!}G_{x}^{-1} = G_{x}^{\dag}$.
For this selection the maps~(\ref{6.10}) form
a Hermitian metric on  $M$.

	Conversely, let there be given an arbitrary Hermitian metric
$g$ with
$G(x) = D^\dag(x)G_{p,q}D(x)$, $D^\dag = D$. Take some constant unitary
matrix $A$ ($A^\dag=A^{-1}$) and any $C$ for which
 $[C]^\dag[C] = A^\dag G_{p,q} A$.
Let us define
$[{^*\!}C] := \left( [C]^\dag \right) ^{-1}$.
Putting $[G_x]=A^\dag D(x)$, from $D^\dag=D$
we get $[G_x]^\dag=[G_x]^{-1}$.
At last, define in a fixed basis
$[{^*\!}G_x] = [G_x]$ (${^*\!}G_{x}^{-1} = G_{x}^{\dag}$).
Thus we have constructed a transport along
$id_{ \{ (1,0),(0,1) \} \times M}$
with  $[F_x]=[C]A^\dag D(x)$ and
$[{^*\!}F_{x}] = \left([C]^\dag\right)^{-1} A^\dag D(x)$.
In particular,
we have $[L_{x\to y}^{1,0}] = D^\dag(y)G_{p,q}D(x)$, so that
 $[L_{x\to x}^{1,0}] =G(x)$.

	In this way we have proved

	\begin{Th}	\label{Th6.1}
	The class of Hermitian metrics on a differentiable manifold is
in one-to-one correspondence with the class of transports along the
identity map in the tensor bundle of rank 1 over that manifold which
have decompositions like (\ref{6.6})--(\ref{6.9})
in which ${^*\!}C^{-1}$ and ${^*\!}G_{x}^{-1}$
(and, hence, ${^*\!}F^{-1}$) are Hermitian conjugate to $C$ and $G_x$
(and $F$) respectively.
	\end{Th}

	\vspace{4ex}\noindent\hspace{4em}%
{\bfseries\large\sffamily\upshape(4)}
	The problems concerning linear transports along maps in the
tensor bundles over a differentiable manifold can be investigated in
the same way as in~\cite[Sect.3]{f-LTP-appl}
in which the text before proposition~3.3 remains true
\emph{mutatis mutandis}
in the considered in the present paper general case.

	\vspace{4ex}\noindent\hspace{4em}%
{\bfseries\large\sffamily\upshape(5)}
	At the end we want to pay attention to tensor densities.
Usually~\cite{Schouten/physics}, a tensor density (field) is defined as
a quantity which is locally represented by a set of numbers (resp.
functions) with a suitable transformation law. Our equivalent view is
that the tensor densities (density fields) are tensors (resp. tensor
fields) that appropriately depend on one fixed basis in the
corresponding tensor space and
which are referred to modified (with respect to the tensors) bases.

	Let $M$ be a differentiable manifold and
a basis $\{ {_0\!}E_{A}^{B}(x) \}$
be fixed in $\left.T_{q}^{p}\right|_x(M)$,
and
$\{ {_1\!}E_{A}^{B}(x) \}$ be an arbitrary basis in it. Here $A$ and
$B$ stand for the corresponding multiindeces (e.g.
$A=(\alpha_1,\ldots,\alpha_p)$, $B=(\beta_1,\ldots,\beta_p)$,
 $\alpha_1,\ldots,\beta_p=1,\ldots,\dim M)$.
We define a tensor density (field) of type $(p,q)$ and weight
$w\in\mathbb{R}$ (with respect to $\{ {_0\!}E_{A}^{B}(x) \}$)
as a tensor (field)
 ${_{0}^{w}\!}\mathcal{T}(x)\in \left.T_{q}^{p}\right|_x(M)$
whose local components ${_{0}^{w}  }\mathcal{T}_{B}^{A}(x)$
are referred to bases like
 $\{ |{_{0}^{1}\!}E(x)|^w {_1\!}E_{A}^{B}(x) \}$, where
 $|{_{0}^{1}\!}E(x)|$ is the Jacobian between the above bases, i.e.
$|{_{0}^{1}\!}E(x)| := \det \left({_{0}^{1}\!}E(x)\right)$,
\(
{_{0}^{1}\!}E(x) :=
\left[\partial x_{1}^{i} /\partial x_{0}^{j}\right].
\)
 So, we have
	\begin{equation}	\label{6.11}
{_{0}^{w}  }\mathcal{T}(x) = {_{0}^{w}  }\mathcal{T}_{B}^{A}(x)
 |{_{0}^{1}\!}E(x)|^w {_1\!}E_{A}^{B}.
	\end{equation}

	It is easy to verify that the components of the so defined
tensor densities have the accepted transformation law~\cite[ch.~II,
Sect.~8]{Schouten/physics}.  Consequently, the both definitions are
equivalent.

	We shall mention only two features of the tensor-density case.

	(i) There exists a class of transports along maps like
 $\eta:\mathbb{R}\times N\to M$, i.e.
 $K_{(v,l)\to(w,m)}^{\eta}$, $v,w\in\mathbb{R},\ l,m\in N$,
which map tensor densities of weight $v$ at one point into such of
weight $w$ at another point. For these transports the results of
Sect.~\ref{III} are valid, in particular,
for $N=M$ and $\eta=id_{\mathbb{R}\times M}$
we have the case considered in
example~\ref{Exam3.1} (with $A=\mathbb{R}$).

	(ii) Of course, one can differentiate a tensor-density field as
tensor field using~(\ref{6.11}), but this operation does not lead
directly to what one expects. In fact, applying~(\ref{4.5}), one finds
	\begin{eqnarray} \nonumber \!\! \!\! \!\! \!\! \!\!
\mathcal{D}_{l^a}^{\varkappa}\left( {_{0}^{w}  }\mathcal{T}  \right)
\!\!&=&\!\!
\left[
\frac{\partial}{\partial l^a}
\left( {_{0}^{w}  }\mathcal{T}_{D}^{C}(\varkappa(l)) \right) +
{_a  }\Gamma_{D\,B}^{C\,A}(l;\varkappa)
\left( {_{0}^{w}  }\mathcal{T}_{A}^{B}(\varkappa(l)) \right)
\right]
\left| {_{0}^{1}\!}E(\varkappa(l)) \right|^w \times
\\	 \label{6.12}
\!\!&\times&\!\! {_1\!}E_{C}^{D}(\varkappa(l)) +
{_{0}^{w}  }\mathcal{T}(\varkappa(l))
\left(
\left| {_{0}^{1}\!}E(\varkappa(l)) \right|^{-w}
\frac{\partial}{\partial l^a}
		\left| {_{0}^{1}\!}E(\varkappa(l)) \right|^w
\right) .
		\end{eqnarray}
The term in parentheses in the last term is equal to
\(
w\frac{\partial}{\partial l^a} \ln
		\left| {_{0}^{1}\!}E(\varkappa(l)) \right|
\)
which, due to~(\ref{4.8}), can be written as
\[
- w \left({_a^1  }\Gamma_{\cdot\,\alpha}^{\alpha\,\cdot}(l;\varkappa) -
	  {_a^0 }\Gamma_{\cdot\,\alpha}^{\alpha\,\cdot}(l;\varkappa)
\right) =
 + w \left({_a^1  }\Gamma_{\alpha\,\cdot}^{\cdot\,\alpha}(l;\varkappa) -
          {_a^0  }\Gamma_{\alpha\,\cdot}^{\cdot\,\alpha}(l;\varkappa),
\right),
\]
where the points in the gamma's stand instead of the absent indices.

Here and above ${_a^1  }\Gamma_{D\,B}^{C\,A}(l;\varkappa)$
are the components of $\mathcal{D}$ in the basis $\{{_1\!}E_{A}^{B}\}$,
i.e.
\[
\mathcal{D}_{l^a}^{\varkappa}(E_{B}^{A}) =
{_a^1  }\Gamma_{D\,B}^{C\,A}(l;\varkappa)
{_1\!}E_{C}^{D}(\varkappa(l)),
\]
and
\(
{_{a}^{0\!}}\Gamma_{\ldots}^{\ldots}:=
{_{a}^{0\!}}\Gamma_{\ldots}^{\ldots}
\left.\right|_{{_1\!}E_{B}^{A} = _0\! E_{B}^{A}}.
\)
	So, if we put
\(
P_{a}^{-}(l;\varkappa):=
{_a^1  }\Gamma_{\cdot \, \alpha}^{\alpha \, \cdot}(l;\varkappa),
\)
\(
P_{a}^{+}(l;\varkappa):=
{_a^1  }\Gamma_{\alpha\,\cdot}^{\cdot\,\alpha}(l;\varkappa),
\)
and
\(
{^0\!}P_{a}^{\pm}(l;\varkappa):=\left.P_{a}^{\pm}(l;\varkappa)
\right|_{_1\!E_{B}^{A} = _0\! E_{B}^{A}},
\)
we get
	\begin{eqnarray} \nonumber
& & \!\! \!\! \!\! \!\! \!\! \!\! \!\!
\mathcal{D}_{l^a}^{\varkappa} ( {_{0}^{w}  }\mathcal{T} )  \pm \>
w \left( {^0\!}P_{a}^{\pm}(l;\varkappa) \right)
\left( {_{0}^{w}  }\mathcal{T}(\varkappa(l)) \right) =
\Bigl[
\frac{\partial}{\partial l^a}
\left({_{0}^{w}  }\mathcal{T}_{D}^{C}(\varkappa(l))\right) +
{_a^1  }\Gamma_{D\,B}^{C\,A}(l;\varkappa) \times
\\ & \times &   \label{6.13}
\left( {_{0}^{w}  }\mathcal{T}_{A}^{B}(\varkappa(l)) \right) \pm
 w P_{a}^{\pm}(l;\varkappa)
\left({_{0}^{w}  }\mathcal{T}_{D}^{C}(\varkappa(l))\right)
 \Bigr]
\left| {_{0}^{1}\!}E(\varkappa(l)) \right|^w
{_1\!}E_{C}^{D}(\varkappa(l)).
	\end{eqnarray}

	Thus the operator~(\ref{4.3}) when applied on tensor density
fields produces, of course, tensor fields which, generally, are not
tensor density fields as by~(\ref{6.12})  their components with respect
to the corresponding bases depend on them in a way different from that
of tensor densities. On the contrary, the right-hand-side
of~(\ref{6.13}) is a tensor density field
whose components, following~\cite[ch.~V, Sect.~1]{Schouten/physics},
should
be identified with those of the \emph{$a$-th partial (plus or minus)
derivation along the map $\varkappa$ } of the initial tensor density
field; the components of this derivation being defined by the r.h.s.
of~(\ref{6.13}).

	It can be proved that when  $\varkappa$ is a path and the
transport along it is a parallel transport assigned to a covariant
differentiation (linear connection) (see~\cite[p.~19]{f-LTP-general}),
the components of the r.h.s. of~(\ref{6.13}) coincide with the covariant
differentiation along the
tangent to the path vector field of the initial tensor density field
(see~\cite[ch.~V, Sect.~1]{Schouten/physics}).

	In the special case when
$\Gamma_{\cdot\,\alpha}^{\alpha\,\cdot} =
\Gamma_{\alpha\,\cdot}^{\cdot\,\alpha}$,
we have $P^-=-P^+$, i.e.
the defined derivation is unique.

	The appropriate approach to the derivation of tensor density
fields is based on transports of tensor densities mentioned in (i) and
the general theory of Sect.~\ref{IV}, but this will be done elsewhere.

\section * {\bfseries Acknowledgement}
\label{VII}
\setcounter {equation} {0}

This work was partially supported by the National Science Foundation
of Bulgaria under Grant No. F642.

\bibliography{bozhopub,bozhoref}

\begin{thebibliography}{10}

\bibitem{f-TP-general}
Bozhidar~Z. Iliev.
\newblock Transports along paths in fibre bundles. {General} theory.
\newblock JINR Communication E5-93-299, Dubna, 1993.

\bibitem{f-LTP-general}
Bozhidar~Z. Iliev.
\newblock Linear transports along paths in vector bundles. {I}.~{General}
  theory.
\newblock JINR Communication E5-93-239, Dubna, 1993.

\bibitem{f-TP-morphisms}
Bozhidar~Z. Iliev.
\newblock Transports along paths in fibre bundles. {III}.~{Consistency} with
  bundle morphisms.
\newblock JINR Communication E5-94-41, Dubna, 1994.
\newblock (LANL xxx archive server, E-print No. dg-ga/9704004).

\bibitem{f-LTP-Cur+Tor-prop}
Bozhidar~Z. Iliev.
\newblock Linear transports along paths in vector bundles. {V}. {Properties} of
  curvature and torsion.
\newblock JINR Communication E5-97-1, Dubna, 1997.
\newblock (LANL xxx archive server, E-print No. dg-ga/9709017).

\bibitem{f-Jacobi}
Bozhidar~Z. Iliev.
\newblock Some generalizations of the {Jacobi} identity with applications to
  the curvature- and torsion-depending hamiltonians of particle systems.
\newblock In J.~{\L}awrynowicz, editor, {\em Hurwitz-type structures and
  applications to surface physics}, number~{II} in Deformations of Mathematical
  Structures, pages 161--188. Kluwer Academic Publishers Group,
  Dordrecht-Boston-London, 1993.
\newblock (Papers from the Seminars on Deformations 1988--1992).

\bibitem{f-LTP-appl}
Bozhidar~Z. Iliev.
\newblock Linear transports along paths in vector bundles. {II}.~{Some}
  applications.
\newblock JINR Communication E5-93-260, Dubna, 1993.

\bibitem{f-LTP-Cur+Tor}
Bozhidar~Z. Iliev.
\newblock Linear transports along paths in vector bundles. {III}.~{Curvature}
  and torsion.
\newblock JINR Communication E5-93-261, Dubna, 1993.

\bibitem{f-LTP-metrics}
Bozhidar~Z. Iliev.
\newblock Linear transports along paths in vector bundles. {IV}.~{Consistency}
  with bundle metrics.
\newblock JINR Communication E5-94-17, Dubna, 1994.

\bibitem{Husemoller}
D.~Husemoller.
\newblock {\em Fibre bundles}.
\newblock McGraw-Hill Book Co., New York-St. Louis-San
  Francisco-Toronto-London-Sydney, 1966.

\bibitem{K&N-1}
S.~Kobayashi and K.~Nomizu.
\newblock {\em Foundations of Differential Geometry}, volume~I.
\newblock Interscience Publishers, New York-London, 1963.

\bibitem{Postnikov-Morse}
M.~M. Postnikov.
\newblock {\em Introduction to the Morse theory}.
\newblock Nauka, Moscow, 1971.
\newblock (In Russian).

\bibitem{Greub&et.al.-1}
W.~Greub, S.~Halperin, and R.~Vanstone.
\newblock {\em De Rham cohomology of manifolds and vector bundles}, volume~1 of
  {\em Connections, Curvature, and Cohomology}.
\newblock Academic Press, New York and London, 1972.

\bibitem{Bellman}
R.~Bellman.
\newblock {\em Introduction to matrix analysis}.
\newblock McGraw-Hill book comp., New York-Toronto-London, 1960.

\bibitem{Schouten/physics}
J.~A. Schouten.
\newblock {\em Tensor analysis for physicists}.
\newblock Clarendon Press, Oxford, 1951.

\end{thebibliography}
\bibliographystyle{unsrt}

\end{document}